\documentclass{article}%
\usepackage{bbm}
\usepackage{epsfig}
\usepackage{graphics,graphicx,amssymb,amsmath,verbatim}
\usepackage{amssymb}
\usepackage{mathrsfs}
\usepackage{amsfonts}
\usepackage{amsmath}
\usepackage{graphicx}
\usepackage{easybmat}
\usepackage{color}%
\providecommand{\U}[1]{\protect \rule{.1in}{.1in}}
\newtheorem{theorem}{Theorem}

\newtheorem{corollary}{Corollary}

\newtheorem{definition}{Definition}

\newtheorem{lemma}{Lemma}

\allowdisplaybreaks[4]

\begin{document}

\title{Improved Razumikhin and Krasovskii Stability Criteria for Time-Varying
Stochastic Time-Delay Systems}
\author{Bin Zhou\thanks{Center for Control Theory and Guidance Technology, Harbin
Institute of Technology, Harbin, 150001, China. Email: binzhoulee@163.com,
binzhou@hit.edu.cn. Corresponding author of this paper.}\quad \quad Weiwei Luo}
\date{}
\maketitle

\begin{abstract}
The problem of $p$th moment stability for time-varying stochastic time-delay
systems with Markovian switching is investigated in this paper. Some novel
stability criteria are obtained by applying the generalized Razumikhin and
Krasovskii stability theorems. Both $p$th moment asymptotic stability and
(integral) input-to-state stability are considered based on the notion and
properties of uniformly stable functions and the improved comparison
principles. The established results show that time-derivatives of the
constructed Razumikhin functions and Krasovskii functionals are allowed to be
indefinite, which improve the existing results on this topic. By applying the
obtained results for stochastic systems, we also analyze briefly the stability
of time-varying deterministic time-delay systems. Finally, examples are
provided to illustrate the effectiveness of the proposed results.
\vspace{0.3cm}


\end{abstract}

\section{Introduction}

Practical control systems need to meet several important features, for
example, time delays, stochastic perturbations, and time-varying parameters.
Time-delay systems refer to a class of dynamics systems whose rate of the
current states is affected by their past states. In most cases, time delays
are negative for the analysis and design of control systems since they may be
the source of performance degradation and instability. On the other hand,
stochastic perturbations are very common in practical problems since many
engineering systems are subject to external random fluctuations including
environmental noise and Markovian switching in the system parameters. Finally,
time-varying control systems are of great importance since in most cases the
system parameters are changing with time and do not fulfill the usual
stationary assumption. There are plenty of papers that deal with control
systems having at least one of these three features. For example, stability
analysis of control systems with time-varying coefficients and/or time-delays
were considered in \cite{czz16auto, Hale1977, Mazenc2015CDC} and
\cite{Zhou2014}; stability analysis of delay-free time-varying stochastic
systems was investigated in \cite{LiuShen2012} and \cite{Wu2006}; stability
analysis of stochastic time-delay systems was studied in \cite{Cong16nahs,
Mao2002tac} and \cite{wssc13tc}; and stability of stochastic neutral
time-delay systems was discussed in \cite{LiZhou2015IJRNC} and
\cite{SPWZ13JFI}.

There are only a few papers that deal with control problems, especially,
stability and stabilization, for systems having all of these three features
mentioned in the above. The Krasovskii functional approach is effective for
handling this class of systems. For example, by this approach, stability
analysis, $H_{\infty}$ control, and $H_{\infty}$ filtering for stochastic
systems with time-delays were respectively investigated in
\cite{YueHan2005tac}, \cite{xlc04scl} and \cite{xlm07tcasi}, and the results
were expressed by linear matrix inequalities. Compared with the Krasovskii
functional approach, the study of the Razumikhin approach was relatively
retarded. Razumikhin-type stability theorems for $p$th moment exponential
stability of time-varying stochastic time-delay systems were originally
established in \cite{Mao1999} and were then utilized to study hybrid
stochastic interval systems in \cite{MaoLam2006}. Later on, Razumikhin-type
theorems for general $p$th moment input-to-state stability and asymptotic
stability of time-varying stochastic time-delay systems were respectively
established in \cite{HuangMao2009TAC} and \cite{HuangMao2008tcsi}. Some other
recently published work related to this issue can be found in
\cite{KangLiu2014tac, Wu2016} and the references therein.

However, in all the mentioned references in the above, as commonly required in
the classical Lyapunov stability theory, time-derivatives of both the
Krasovskii functional and the Razumikhin function need to be negative definite
(under the Razumikhin condition for the later case), which may be
conservative. Such a requirement was relaxed in \cite{Wu2006} where a
time-varying function was introduced so that the time-derivative of the
Lyapunov function along the trajectories of a delay-free stochastic system can
be indefinite. This idea was also utilized in \cite{pz10tac} and
\cite{pz10spl} to build Razumikhin-type stability theorems for time-varying
stochastic time-delay systems. This method, together with some other new
ideas, has also been utilized in \cite{ning13scl} and \cite{Ning2014tac,
Ning2015AMC} to study stability of deterministic time-varying nonlinear
systems and time-delay systems, respectively. A different method that also
allows indefinite time-derivatives of the Lyapunov function was established
recently in \cite{Wu2013}.

Very recently, to weaken the condition required by the classical Lyapunov
stability theory that the time-derivative of the Lyapunov function along the
trajectories of a delay-free linear time-varying system is negative definite,
a new Lyapunov approach that allows the time-derivative of the Lyapunov
function to take both negative and positive values was proposed in
\cite{Zhou16AUTO} by using the notion of scalar stable scalar functions.
Motivated by the results of \cite{Zhou16AUTO}, the classical Razumikhin and
Krasovskii stability theorems were generalized in \cite{ZE16AUTO} for
time-varying time-delay systems by allowing time-derivatives of the Razumikhin
functions and the Krasovskii functionals to be indefinite. The approach was
also utilized in \cite{Zhou15arXiv} to study input-to-state stability of
general nonlinear time-varying systems. As mentioned in \cite{Zhou15arXiv} and
\cite{ZE16AUTO}, the main advantage of the approach in \cite{Zhou15arXiv} and
\cite{ZE16AUTO} is that the scalar function there does not need to satisfy
some restrictive conditions that were required in
\cite{ning13scl,Ning2014tac,Ning2015AMC,pz10tac,Wu2006} and \cite{Wu2013} (see
detailed explanations at the end of Section \ref{sec3}).

Motivated by these existing work, especially our recent work \cite{Zhou16AUTO,
Zhou15arXiv} and \cite{ZE16AUTO}, in the present paper we will provide some
new Razumikhin-type and Krasovskii-type stability theorems for time-varying
stochastic time-delay systems. Both $p$th moment (integral) input-to-state
stability and asymptotic stability will be considered. The established results
possess the significant feature that the time-derivatives of the Razumikhin
functions and Krasovskii functionals are not required to be negative definite,
which relaxes the existing results on the same problems. To derive our
results, the comparison lemmas build in \cite{Zhou15arXiv} and \cite{ZE16AUTO}
were generalized to the case that the right hand side of the comparison
differential equation contains a drift term. With this new comparison lemma,
stability theorems are established based on the generalized It\^{o} formula
and stochastic analysis theory. The improvement of the proposed results over
the existing ones, especially those provided in \cite{pz10tac} and
\cite{Wu2013} will be made clear. Applications of the established theorems to
some special systems such as time-varying stochastic time-delay systems
without Markovian jumping and time-varying deterministic time-delay systems
are also discussed in detail. We also provide some numerical examples to
illustrate the effectiveness of the proposed theorems.

The remaining of the paper is organized as follows. Problem formulation and
some preliminary results will be presented in Section \ref{sec2}.
Razumikhin-type and Krasovskii-type stability theorems for time-varying
stochastic time-delay systems are respectively presented in Sections
\ref{sec3} and \ref{sec4}. Applications of the theory established in Sections
\ref{sec3} and \ref{sec4} to time-varying deterministic time-delay systems
will be shown in Section \ref{sec5}. Numerical examples are given in Section
\ref{sec6} and finally Section \ref{sec7} concludes the paper.

\section{\label{sec2}Problem Formulation and Preliminaries}

\subsection{Problem Formulation}

Throughout this paper, unless otherwise specified, $\mathbf{R}^{n}$ and
$\mathbf{R}^{n\times m}$ denote, respectively, the $n$-dimensional Euclidean
space and the set of $n\times m$ real matrices. Let $(\Omega,\mathcal{F}%
,\left \{  \mathcal{F}\right \}  _{t\geq t_{0}},\mathbb{P})$ be a complete
probability space with a filtration $\left \{  \mathcal{F}\right \}  _{t\geq
t_{0}}$ satisfying the usual conditions (\emph{i.e.}, it is right continuous
and $\mathcal{F}_{t_{0}}$ contains all $\mathbb{P}$-null sets) and
$\mathbf{E}\left \{  \cdot \right \}  $ be the expectation operator with respect
to the probability measure. Let $w\left(  t\right)  =\left[  w_{1}\left(
t\right)  ,\cdots,w_{m}\left(  t\right)  \right]  ^{\mathrm{T}}$ be an
$m$-dimensional Brownian motion defined on the probability space. Let
$J=[0,\infty)$, $\left \vert \cdot \right \vert $ denote the Euclidean norm in
$\mathbf{R}^{n}$, and $\left \Vert f\right \Vert _{[t_{0},t]}=\sup \left \{
\left \vert f\left(  s\right)  \right \vert ,s\in \lbrack t_{0},t]\subset
J\right \}  .$ Let $\mathbf{C}\left(  \left[  -\tau,0\right]  ,\mathbf{R}%
^{n}\right)  $, where $\tau \geq0$, denote the family of all continuous
$\mathbf{R}^{n}$-valued function $\varepsilon$ defined on $[-\tau,0]$ with the
norm $\left \Vert \epsilon \right \Vert =\sup_{-\tau \leq \theta \leq0}\left \vert
\epsilon \left(  \theta \right)  \right \vert .$ Let $\mathbf{C}_{\mathcal{F}%
_{t_{0}}}^{\mathrm{b}}\left(  \left[  -\tau,0\right]  ,\mathbf{R}^{n}\right)
$ be the family of all $\mathcal{F}_{t_{0}}$-measurable bounded $\mathbf{C}%
\left(  \left[  -\tau,0\right]  ,\mathbf{R}^{n}\right)  $-valued random
variables $\xi=\left \{  \xi \left(  \theta \right)  :-\tau \leq \theta
\leq0\right \}  .$ For $p>0$ and $t\geq t_{0},$ denote by $L_{\mathcal{F}_{t}%
}^{p}\left(  \left[  -\tau,0\right]  ,\mathbf{R}^{n}\right)  $ the family of
all $\mathcal{F}_{t}$-measurable $\mathbf{C}\left(  \left[  -\tau,0\right]
,\mathbf{R}^{n}\right)  $-valued random processes $\eta=\left \{  \eta \left(
\theta \right)  :-\tau \leq \theta \leq0\right \}  $ such that $\sup_{-\tau
\leq \theta \leq0}\mathbf{E}\left \{  \left \vert \eta \left(  \theta \right)
\right \vert ^{p}\right \}  <\infty.$ We denote by $\lambda \in \mathcal{VK}%
_{\infty}$ and $\lambda \in \mathcal{CK}_{\infty}$ if $\lambda \in \mathcal{K}%
_{\infty}$ and $\lambda$ is convex and concave, respectively. We use
$\mathbf{C}\left(  J,\mathbf{R}^{n}\right)  $ and $\mathbf{PC}\left(
J,\mathbf{R}^{n}\right)  $ to denote respectively the space of $\mathbf{R}%
^{n}$-valued continuous functions and piecewise continuous functions defined
on $J,$ and $J_{\tau}=[-\tau,\infty).$ We also denote $\mathbf{L}_{\infty}%
^{l}\left(  J\right)  =\left \{  f\left(  t\right)  :J\rightarrow \mathbf{R}%
^{l},|\sup_{t\in J}\left \{  \left \vert f\left(  t\right)  \right \vert
\right \}  <\infty \right \}  .$

Let $r\left(  t\right)  ,t\in J$ be a right-continuous Markov chain on the
probability space taking values in a finite state $\mathbf{S}=\left \{
1,2,\cdots,N\right \}  $ with generator $\Gamma=\left(  \gamma_{ij}\right)
_{N\times N}$ given by
\[
\mathbf{P}\left \{  r\left(  t+\Delta \right)  =j|r\left(  t\right)  =i\right \}
=\left \{
\begin{array}
[c]{r}%
\gamma_{ij}\Delta+o\left(  \Delta \right) , \text{ \  \ if }i\neq j,\\
1+\gamma_{ij}\Delta+o\left(  \Delta \right) , \text{ \ if }i=j,
\end{array}
\right.
\]
where $\Delta>0$ and $\gamma_{ij}\geq0$ is the transition rate from $i$ to $j$
if $i\neq j$ while $\gamma_{ii}=-\sum_{i\neq j}\gamma_{ij}.$ Assume that the
Markov chain $r\left(  \cdot \right)  $ is independent of the Brownian motion
$w\left(  \cdot \right)  $. It is known that almost all sample paths of
$r\left(  t\right)  $ are right-continuous step functions with a finite number
of simple jumps in any finite subinterval of $J$.

We consider the following time-varying stochastic functional differential
equation (SFDE) with Markovian switching \cite{HuangMao2009TAC}:%
\begin{equation}
\mathrm{d}x\left(  t\right)  =f\left(  t,x_{t},r\left(  t\right)  ,u\left(
t\right)  \right)  \mathrm{d}t+g\left(  t,x_{t},r\left(  t\right)  ,u\left(
t\right)  \right)  \mathrm{d}w\left(  t\right)  ,\text{ }t\in J, \label{sys}%
\end{equation}
where the initial state is $x_{t_{0}}=\xi \in \mathbf{C}_{\mathcal{F}_{t_{0}}%
}^{b}\left(  \left[  -\tau,0\right]  ,\mathbf{R}^{n}\right)  $, the input
$u:J\rightarrow \mathbf{R}^{l}$ is assumed to be locally essentially bounded,
and $x_{t}=x\left(  t+\theta \right)  ,\theta \in \left[  -\tau,0\right]  $ is
regarded as a $\mathbf{C}$-valued stochastic process. We assume that
\begin{align*}
f  &  :J\times \mathbf{C}\left(  \left[  -\tau,0\right]  ,\mathbf{R}%
^{n}\right)  \times \mathbf{S}\times \mathbf{R}^{l}\rightarrow \mathbf{R}^{n},\\
g  &  :J\times \mathbf{C}\left(  \left[  -\tau,0\right]  ,\mathbf{R}%
^{n}\right)  \times \mathbf{S}\times \mathbf{R}^{l}\rightarrow \mathbf{R}%
^{n\times m},
\end{align*}
are measurable functions with $f\left(  t,0,i,0\right)  \equiv0$ and $g\left(
t,0,i,0\right)  \equiv0$ for all $t\in J$ and $i\in \mathbf{S}$, and are
sufficiently smooth so that (\ref{sys}) only has continuous solution on $J$.
Thus, (\ref{sys}) admits a trivial solution $x\left(  t,t_{0},0\right)
\equiv0,t\geq t_{0}\in J.$ In the absence of Markovian jumping, the SFDE
(\ref{sys}) degrades into the following SFDE \cite{HuangDeng2009}%
\begin{equation}
\mathrm{d}x\left(  t\right)  =f_{0}\left(  t,x_{t},u(t)\right)  \mathrm{d}%
t+g_{0}\left(  t,x_{t},u(t)\right)  \mathrm{d}w\left(  t\right)  ,\text{ }t\in
J, \label{sys2}%
\end{equation}
where $f_{0}:J\times \mathbf{C}\left(  \left[  -\tau,0\right]  ,\mathbf{R}%
^{n}\right)  \rightarrow \mathbf{R}^{n}$ and $g_{0}:J\times \mathbf{C}\left(
\left[  -\tau,0\right]  ,\mathbf{R}^{n}\right)  \rightarrow \mathbf{R}^{n\times
m}$ are measurable functions with $f_{0}\left(  t,0,0\right)  \equiv0$ and
$g_{0}\left(  t,0,0\right)  \equiv0$ for all $t\in J,$ and are sufficiently
smooth so that (\ref{sys2}) only has continuous solution on $J$.

In this paper we are interested in the stability analysis of the SFDE
(\ref{sys}) and (\ref{sys2}). To this end, we give the following definitions.

\begin{definition}
Let $p>0$ be a constant. The trivial solution of the SFDE (\ref{sys}) or
(\ref{sys2}) is said to be

\begin{enumerate}
\item $p$th moment input-to-state stable (ISS) if there exist $\sigma
\in \mathcal{KL}$ and $\gamma_{1}\in \mathcal{K}$ such that, for any
$u\in \mathbf{L}_{\infty}^{l},$
\begin{equation}
\mathbf{E}\left \{  \left \vert x\left(  t\right)  \right \vert ^{p}\right \}
\leq \sigma \left(  \mathbf{E}\left \{  \left \Vert \xi \right \Vert ^{p}\right \}
,t-t_{0}\right)  +\gamma_{1}\left(  \left \Vert u\right \Vert _{[t_{0}%
,t]}\right)  ,\;t,t_{0}\in J,\;t\geq t_{0}. \label{pISS}%
\end{equation}

\item $p$th moment integral input-to-state stable (iISS) if there exist
$\sigma \in \mathcal{KL}$ and $\gamma_{1},\gamma_{2}\in \mathcal{K}$ such that,
for any $u,$%
\[
\mathbf{E}\left \{  \left \vert x\left(  t\right)  \right \vert ^{p}\right \}
\leq \sigma \left(  \mathbf{E}\left \{  \left \Vert \xi \right \Vert ^{p}\right \}
,t-t_{0}\right)  +\gamma_{1}\left(  \int_{t_{0}}^{t}\gamma_{2}\left(
\left \vert u\left(  s\right)  \right \vert \right)  \mathrm{d}s\right)  ,\text{
}t,t_{0}\in J,\;t\geq t_{0}.
\]

\end{enumerate}
\end{definition}

\begin{definition}
Let $p>0$ be a constant. The trivial solution of the SFDE (\ref{sys}) or
(\ref{sys2}) with $u(t)\equiv0$ is said to be

\begin{enumerate}
\item $p$th moment globally uniformly asymptotically stable (GUAS), if there
exists a $\mathcal{KL}$-function $\sigma$ such that%
\[
\mathbf{E}\left \{  \left \vert x\left(  t\right)  \right \vert ^{p}\right \}
\leq \sigma \left(  \mathbf{E}\left \{  \left \Vert \xi \right \Vert ^{p}\right \}
,t-t_{0}\right)  ,\;t,t_{0}\in J,\;t\geq t_{0}.
\]

\item $p$th moment globally uniformly exponentially stable (GUES), if there
exist two constants $\alpha>0$ and $\beta>0$ such that%
\[
\mathbf{E}\left \{  \left \vert x\left(  t\right)  \right \vert ^{p}\right \}
\leq \beta \mathrm{e}^{-\alpha \left(  t-t_{0}\right)  }\mathbf{E}\left \{
\left \Vert \xi \right \Vert ^{p}\right \}  ,\;t,t_{0}\in J,\;t\geq t_{0}.
\]

\end{enumerate}
\end{definition}

The concepts of ISS and iISS, originally introduced in \cite{Sontag1989tac},
have received much attention due to their wide usages in characterizing the
effects of external inputs on a control system. The (i)ISS implies that, no
matter what the size of the initial state is, the state will eventually
approach to a neighborhood of the origin whose size is proportional to the
magnitude of the input. The ISS property is frequently characterized by the
ISS-Lyapunov function \cite{Sontag1998}. As usual, the time-derivative of the
ISS-Lyapunov function is required to be negative definite under some
additional condition on the input signal (\textit{e.g.},
\cite{HuangMao2009TAC, Wu2016}). In this paper, we will weaken such a condition.

At the end of this section, we introduce the definition of infinitesimal
operator for a Lyapunov functional. For a given functional $V:J\times
\mathbf{C}\left(  \left[  -\tau,0\right]  ,\mathbf{R}^{n}\right)
\times \mathbf{S\rightarrow R,}$ its infinitesimal operator $\mathcal{L}$ is
defined by \cite{Mao2002tac}
\[
\mathcal{L}_{V}\left(  t,x_{t},i\right)  =\lim_{\Delta \rightarrow0^{+}}%
\frac{1}{\Delta}\left[  \mathbf{E}\left \{  V\left(  t,x_{t+\Delta},r\left(
t+\Delta \right)  \right)  |x_{t},r\left(  t\right)  =i\right \}  -V\left(
t,x_{t},i\right)  \right]  .
\]
From the above definition, if
\[
V\left(  t,\phi,i\right)  =V\left(  t,\phi(0),i\right)  ,\text{ }\phi
\in \mathbf{C}\left(  \left[  -\tau,0\right]  ,\mathbf{R}^{n}\right)
\mathbf{,}%
\]
where $V\left(  t,x,i\right)  \in \mathbf{C}^{2,1}\left(  J\times \mathbf{R}%
^{n}\times \mathbf{S},J\right)  ,$ in which $\mathbf{C}^{2,1}\left(
J\times \mathbf{R}^{n}\times S,J\right)  $ denotes the family of all
nonnegative functions $V\left(  t,x,i\right)  $ on $J\times \mathbf{R}%
^{n}\times \mathbf{S}$ that are twice continuously differentiable in $x$ and
once in $t,$ then the infinitesimal operator of $V\left(  t,x,i\right)  $
along the SFDE (\ref{sys}) is given by \cite{Mao1999}
\begin{align}
\mathcal{L}_{V}\left(  t,x_{t},i\right)   &  =V_{t}\left(  t,x,i\right)
+V_{x}\left(  t,x,i\right)  f\left(  t,x_{t},i,u\right) \nonumber \\
&  +\frac{1}{2}\mathrm{trace}\left(  g^{\mathrm{T}}\left(  t,x_{t},i,u\right)
V_{xx}\left(  t,x,i\right)  g\left(  t,x_{t},i,u\right)  \right)
+\sum \limits_{j=1}^{N}\gamma_{ij}V\left(  t,x,j\right)  , \label{LV}%
\end{align}
where%
\begin{align*}
V_{t}\left(  t,x,i\right)   &  =\frac{\partial V\left(  t,x,i\right)
}{\partial t},\\
V_{x}\left(  t,x,i\right)   &  =\left[  \frac{\partial V\left(  t,x,i\right)
}{\partial x_{1}},\cdots,\frac{\partial V\left(  t,x,i\right)  }{\partial
x_{n}}\right]  ,\\
V_{xx}\left(  t,x,i\right)   &  =\left(  \frac{\partial^{2}V\left(
t,x,i\right)  }{\partial x_{i}\partial x_{j}}\right)  _{n\times n}.
\end{align*}

\subsection{USF and Comparison Principles}

To build our results, we need the following concept of uniformly stable
functions which are recalled from \cite{Zhou16AUTO} and \cite{ZE16AUTO}.

\begin{definition}
\cite{Zhou16AUTO} A function $\mu \in \mathbf{C}\left(  J,\mathbf{R}\right)  $
is said to be a uniformly stable function (USF) if the following linear
time-varying (LTV) system is GUAS:%
\[
\dot{y}\left(  t\right)  =\mu \left(  t\right)  y\left(  t\right)  ,\text{
}t\in J.
\]
Hence, $\mu \left(  t\right)  $ is a USF if and only if there exist two
constants $\varepsilon>0$ and $\delta \geq0$ such that
\begin{equation}
\int_{t_{0}}^{t}\mu \left(  s\right)  \mathrm{d}s\leq-\varepsilon \left(
t-t_{0}\right)  +\delta,\text{ }t,t_{0}\in J,t\geq t_{0}. \label{USF}%
\end{equation}

\end{definition}

The following two concepts are also important in developing our
Razumikhin-type stability theorems.

\begin{definition}
\cite{ZE16AUTO} Let $\mu \left(  t\right)  $ be a USF. Then the set
\[
\Omega_{\mu}=\left \{  T>0:\sup_{t\in J}\left \{  \int_{t}^{t+T}\mu \left(
s\right)  \mathrm{d}s\right \}  <0\right \}  ,
\]
is said to be the uniform convergence set (UCS) of $\mu \left(  t\right)  $.
\end{definition}

\begin{definition}
\cite{ZE16AUTO} Let $\mu \left(  t\right)  $ be a USF. Then for any given
$T\geq0$ the overshoot $\varphi_{\mu}\left(  T\right)  $ of $\mu \left(
t\right)  $ is defined as%
\[
\varphi_{\mu}\left(  T\right)  =\sup_{t\in J}\left \{  \max_{\theta \in
\lbrack0,T]}\left \{  \int_{t}^{t+\theta}\mu \left(  s\right)  \mathrm{d}%
s\right \}  \right \}  .
\]

\end{definition}

From Proposition 2 in \cite{ZE16AUTO}, we can see that $\varphi_{\mu}\left(
T\right)  $ is a non-decreasing function of $T$ and $0\leq \varphi_{\mu}\left(
T\right)  \leq \delta<\infty$, where $\delta$ is from (\ref{USF}).

Next, we introduce some improved comparison principles for time-varying systems.

\begin{lemma}
\label{GGBI}\cite{Flett1980book} Assume that $\pi \in \mathbf{PC}\left(
J,\mathbf{R}\right)  ,\mu \in \mathbf{C}\left(  J,\mathbf{R}\right)  $ and
$y:J\rightarrow J$ is such that
\[
\mathbf{D}^{+}y\left(  t\right)  \leq \mu \left(  t\right)  y\left(  t\right)
+\pi \left(  t\right)  ,\text{ }t\in J.
\]
Then, for any $t\geq s\in J,$ the following inequality holds true%
\[
y\left(  t\right)  \leq y\left(  s\right)  \exp \left(  \int_{s}^{t}\mu \left(
r\right)  \mathrm{d}r\right)  +\int_{s}^{t}\exp \left(  \int_{\lambda}^{t}%
\mu \left(  r\right)  \mathrm{d}r\right)  \pi \left(  \lambda \right)
\mathrm{d}\lambda.
\]

\end{lemma}

\begin{lemma}
\label{Comparison}Let $y\left(  t\right)  $ $:J\rightarrow J$ be a continuous
function such that
\begin{equation}
\mathbf{D}^{+}y\left(  t\right)  \leq \mu \left(  t\right)  y\left(  t\right)
+\pi \left(  t\right)  \text{ whenever }y\left(  t\right)  \geq \psi \left(
t\right)  , \label{eq10}%
\end{equation}
where $\mu \in \mathbf{C}\left(  J,\mathbf{R}\right)  $ is a USF, $\pi
\in \mathbf{PC}\left(  J,J\right)  $ and $\psi \left(  \cdot \right)  $ is a
$\mathcal{K}$-function. Then, for any given constant $T>0$ and any $t\geq T,$
there holds%
\[
y\left(  t\right)  \leq \max \left \{  y\left(  t-T\right)  \exp \left(
\int_{t-T}^{t}\mu \left(  s\right)  \mathrm{d}s\right)  ,\sup_{t-T\leq s\leq
t}\left \{  \psi \left(  s\right)  \right \}  \mathrm{e}^{\varphi_{\mu}\left(
T\right)  }\right \}  +\int_{t-T}^{t}\exp \left(  \int_{s}^{t}\mu \left(
r\right)  \mathrm{d}r\right)  \pi \left(  s\right)  \mathrm{d}s.
\]

\end{lemma}

Letting $\pi(t)=0$ and $\mathbf{D}^{+}y\left(  t\right)  =\dot{y}\left(
t\right)  $ in Lemma \ref{Comparison} gives immediately Lemma 4 in
\cite{ZE16AUTO}. These two lemmas in the above will play critical roles in
establishing respectively Krasovskii-type and Razumikhin-type stability
theorems in the next two sections.

\section{\label{sec3}Razumikhin-Type Stability Theorems}

This section is concerned with Razumikhin-type stability theorems for the
time-varying SFDE (\ref{sys}). We first give the following result regarding
ISS of (\ref{sys}).

\begin{theorem}
\label{th1}(Razumikhin Theorem for ISS) The SFDE (\ref{sys}) is $p$th moment
ISS if there exist a $V(t,x,i)\in \mathbf{C}^{2,1}\left(  J_{\tau}%
\times \mathbf{R}^{n}\times \mathbf{S},J\right)  $, a USF $\mu(t)$, $\alpha
_{1}\in \mathcal{VK}_{\infty},\alpha_{2}\in \mathcal{CK}_{\infty},q\in
\mathcal{K}_{\infty},\varpi \in \mathcal{K}_{\infty}$, a number $T\in \Omega
_{\mu}$, and a constant $\rho \in(0,1)$ satisfying
\begin{equation}
q\left(  s\right)  \geq \frac{\mathrm{e}^{\varphi_{\mu}(T)}}{\rho}%
s,\quad \forall s\geq0, \label{qs}%
\end{equation}
such that the following conditions are met for all $t\in J$:
\begin{align}
\mathrm{(A}_{1}\mathrm{).}\; \;  &  \alpha_{1}\left(  |x|^{p}\right)  \leq
V\left(  t,x,i\right)  \leq \alpha_{2}\left(  |x|^{p}\right)  ,\quad \forall
i\in \mathbf{S.}\nonumber \\
\mathrm{(B}_{1}\mathrm{).}\; \;  &  \mathbf{E}\left \{  \mathcal{L}_{V}\left(
t,x_{t},i\right)  \right \}  \leq \mu(t)\mathbf{E}\left \{  V\left(
t,x(t),i\right)  \right \}  ,\; \forall i\in \mathbf{S,}\nonumber \\
&  \quad \text{if}\; \max_{i\in \mathbf{S}}\left \{  \mathbf{E}\left \{  V\left(
t,x(t),i\right)  \right \}  \right \}  \geq \max \left \{  \varpi \left(  \left \vert
u\left(  t\right)  \right \vert \right)  ,q^{-1}\left(  \min_{i\in \mathbf{S}%
}\left \{  \mathbf{E}\left \{  V\left(  t+\theta,x(t+\theta),i\right)  \right \}
\right \}  \right)  \right \}  ,\; \forall \theta \in \lbrack-\tau,0]\mathbf{.}
\label{ISSEV}%
\end{align}

\end{theorem}

We next present the following result regarding the characterization of iISS.

\begin{theorem}
\label{th2}(Razumikhin Theorem for iISS) The SFDE (\ref{sys}) is $p$th moment
iISS if there exist a $V(t,x,i)\in \mathbf{C}^{2,1}\left(  J_{\tau}%
\times \mathbf{R}^{n}\times \mathbf{S},J\right)  $, a USF $\mu(t)$, $\alpha
_{1}\in \mathcal{VK}_{\infty},\alpha_{2}\in \mathcal{CK}_{\infty},\varpi
\in \mathcal{K}_{\infty}$, a number $T\in \Omega_{\mu}$, $q\in \mathcal{K}%
_{\infty}$ and a constant $\rho \in(0,1)$ satisfying (\ref{qs}) such that
$\mathrm{(A}_{1}\mathrm{)}$ and the following conditions are met for all $t\in
J$:
\begin{align*}
\mathrm{(B}_{2}\mathrm{).}\; \;  &  \mathbf{E}\left \{  \mathcal{L}_{V}\left(
t,x_{t},i\right)  \right \}  \leq \mu(t)\mathbf{E}\left \{  V\left(
t,x(t),i\right)  \right \}  +\varpi \left(  \left \vert u\left(  t\right)
\right \vert \right)  ,\; \forall i\in \mathbf{S,}\\
&  \quad \text{if}\; \min_{j\in \mathbf{S}}\left \{  \mathbf{E}\left \{  V\left(
t+\theta,x(t+\theta),j\right)  \right \}  \right \}  \leq q\left(  \max
_{i\in \mathbf{S}}\left \{  \mathbf{E}\left \{  V\left(  t,x(t),i\right)
\right \}  \right \}  \right)  ,\forall \theta \in \left[  -\tau,0\right]
\mathbf{.}%
\end{align*}

\end{theorem}

Setting $u(t)\equiv0$ in Theorems \ref{th1} and \ref{th2} gives immediately
the following corollary regarding Razumikhin-type theorem for the $p$th moment
GUAS of the SFDE (\ref{sys}).

\begin{corollary}
\label{GUAS}(Razumikhin Theorem for GUAS) The SFDE (\ref{sys}) is $p$th moment
GUAS if there exist a $V(t,x,i)\in \mathbf{C}^{2,1}\left(  J_{\tau}%
\times \mathbf{R}^{n}\times \mathbf{S},J\right)  $, a USF $\mu(t)$, $\alpha
_{1}\in \mathcal{VK}_{\infty},\alpha_{2}\in \mathcal{CK}_{\infty},q\in
\mathcal{K}_{\infty}$, a number $T\in \Omega_{\mu}$, and a constant $\rho
\in(0,1)$ satisfying (\ref{qs}) such that $\mathrm{(A}_{1}\mathrm{)}$ and the
following conditions are met for all $t\in J$:
\begin{align}
\mathrm{(B}_{3}\mathrm{)}  &  .\; \mathbf{E}\left \{  \mathcal{L}_{V}\left(
t,x_{t},i\right)  \right \}  \leq \mu(t)\mathbf{E}\left \{  V\left(
t,x(t),i\right)  \right \}  ,\; \forall i\in \mathbf{S},\nonumber \\
&  \quad \text{if}\; \min_{j\in \mathbf{S}}\left \{  \mathbf{E}\left \{  V\left(
t+\theta,x(t+\theta),j\right)  \right \}  \right \}  \leq q\left(  \max
_{i\in \mathbf{S}}\left \{  \mathbf{E}\left \{  V\left(  t,x(t),i\right)
\right \}  \right \}  \right)  ,\; \forall \theta \in \left[  -\tau,0\right]
\mathbf{.} \label{EVtq1}%
\end{align}
The SFDE is furthermore $p$th moment GUES if there exist three positive
numbers $\beta_{i},i=0,1,2$ such that the following condition is satisfied:%
\[
\mathrm{(D).}\; \alpha_{1}(s)=\beta_{1}s^{\beta_{0}},\; \alpha_{2}%
(s)=\beta_{2}s^{\beta_{0}}.
\]

\end{corollary}

Applying Theorems \ref{th1} and \ref{th2} and Corollary \ref{GUAS} on system
(\ref{sys2}) gives the following corollary.

\begin{corollary}
\label{SFDE}Consider the SFDE (\ref{sys2}). Let $V(t,x)\in \mathbf{C}%
^{2,1}\left(  J_{\tau}\times \mathbf{R}^{n},J\right)  $, $\mu(t)$ be a USF,
$\alpha_{1}\in \mathcal{VK}_{\infty},\alpha_{2}\in \mathcal{CK}_{\infty}%
,\varpi \in \mathcal{K}_{\infty},q\in \mathcal{K}_{\infty}$, $T\in \Omega_{\mu}$,
and $\rho \in(0,1)$ be a constant satisfying (\ref{qs}). Consider the following
conditions, where $t\in J$:
\begin{align*}
\mathrm{(A}_{2}\mathrm{)}  &  .\; \; \alpha_{1}\left(  |x|^{p}\right)  \leq
V\left(  t,x\right)  \leq \alpha_{2}\left(  |x|^{p}\right)  .\\
\mathrm{(B}_{4}\mathrm{)}  &  .\; \; \mathbf{E}\left \{  \mathcal{L}_{V}\left(
t,x_{t}\right)  \right \}  \leq \mu(t)\mathbf{E}\left \{  V\left(  t,x(t)\right)
\right \}  \mathbf{,}\quad \\
&  \text{if}\; \mathbf{E}\left \{  V\left(  t,x(t)\right)  \right \}  \geq
\max \left \{  \varpi \left(  \left \vert u\left(  t\right)  \right \vert \right)
,q^{-1}\left(  \mathbf{E}\left \{  V\left(  t+\theta,x(t+\theta)\right)
\right \}  \right)  \right \}  ,\; \forall \theta \in \lbrack-\tau,0]\mathbf{.}\\
\mathrm{(B}_{5}\mathrm{)}  &  .\; \; \mathbf{E}\left \{  \mathcal{L}_{V}\left(
t,x_{t}\right)  \right \}  \leq \mu(t)\mathbf{E}\left \{  V\left(  t,x(t)\right)
\right \}  +\varpi \left(  \left \vert u\left(  t\right)  \right \vert \right) ,\\
&  \text{if}\; \mathbf{E}\left \{  V\left(  t+\theta,x(t+\theta)\right)
\right \}  \leq q\left(  \mathbf{E}\left \{  V\left(  t,x(t)\right)  \right \}
\right)  ,\; \forall \theta \in \left[  -\tau,0\right]  .\\
\mathrm{(B}_{6}\mathrm{)}  &  .\; \; \mathbf{E}\left \{  \mathcal{L}_{V}\left(
t,x_{t}\right)  \right \}  \leq \mu(t)\mathbf{E}\left \{  V\left(  t,x(t)\right)
\right \}  ,\\
&  \text{if}\; \mathbf{E}\left \{  V\left(  t+\theta,x(t+\theta)\right)
\right \}  \leq q\left(  \mathbf{E}\left \{  V\left(  t,x(t)\right)  \right \}
\right)  ,\; \forall \theta \in \left[  -\tau,0\right]  .
\end{align*}
Then the SFDE (\ref{sys2})

\begin{enumerate}
\item is $p$th moment ISS if $\mathrm{(A}_{2}\mathrm{)}$ and $\mathrm{(B}%
_{4}\mathrm{)}$ are satisfied.

\item is $p$th moment iISS if $\mathrm{(A}_{2}\mathrm{)}$ and $\mathrm{(B}%
_{5}\mathrm{)}$ are satisfied.

\item is $p$th moment GUAS if $\mathrm{(A}_{2}\mathrm{)}$ and $\mathrm{(B}%
_{6}\mathrm{)}$ are satisfied.

\item is $p$th moment GUES if $\mathrm{(A}_{2}\mathrm{)}$, $\mathrm{(B}%
_{6}\mathrm{)}$ and $\mathrm{(D)}$ are satisfied.
\end{enumerate}
\end{corollary}

The classical Razumikhin-type stability theorems for deterministic systems
were generalized to SFDEs in \cite{HuangMao2009TAC, HuangMao2008tcsi,
HuangDeng2009, Mao1999}. In these Razumikhin-stability theorems,
time-derivatives of the Razumikhin functions must be negative definite. Such a
requirement was relaxed in \cite{Wu2006}, where the $p$th moment stability for
delay-free stochastic differential equations was investigated by allowing
time-derivatives of Lyapunov functions to be indefinite. The same idea was
utilized in \cite{pz10tac} to study asymptotic stability of SFDEs and some
improved Razumikhin-type theorems, which take similar forms as Corollary
\ref{GUAS}, were established. However, the corresponding function $\mu \left(
t\right)  $ there needs to satisfy
\begin{equation}
\int_{0}^{\infty}\max \left \{  \mu \left(  t\right)  ,0\right \}  \mathrm{d}%
t<\infty. \label{restrictive}%
\end{equation}
This condition is very restrictive since it cannot be satisfied by any
continuous periodic functions having positive values (see more explanations in
Remark 3 in \cite{ZE16AUTO} and the examples given in Section \ref{sec6}).
Though the restrictive condition (\ref{restrictive}) was not required in
\cite{Wu2013}, another condition%
\begin{equation}
\mu \left(  t\right)  \geq-\frac{\ln q}{\tau},\text{ }t\geq0,
\label{restrictive2}%
\end{equation}
was imposed there, where $q>1$ is some constant and $\tau$ is the delay. We
mention that the restrictive condition (\ref{restrictive}) was also required
in some other related references such as \cite{ning13scl,Ning2014tac} and
\cite{Ning2015AMC}. In the proposed generalized Razumikhin-type stability
theorems (Theorem \ref{th1}, \ref{th2}, and Corollary \ref{GUAS}), both of
these two restrictions (\ref{restrictive})-(\ref{restrictive2}) are not
required. In this sense we emphasize that Corollary \ref{GUAS} generalizes
Theorem 3.2 in \cite{pz10tac} and Item 3 of Corollary \ref{SFDE} generalizes
the results in \cite{Mao2007}. Finally, we mention that, to the best of our
knowledge, the Razumikhin-type stability criteria for (i)ISS of SFDFs without
Markovian switching shown in Corollary \ref{SFDE} was not available in the
literature even if $\mu(t)$ is a negative constant.

\section{\label{sec4}Krasovskii-Type Stability Theorems}

In this section, the Krasovskii functionals with indefinite time-derivatives
are employed to derive $p$th moment ISS, iISS, GUAS and GUES of systems
(\ref{sys}) and (\ref{sys2}).

\begin{theorem}
\label{th3}(Krasovskii Theorem for ISS) The SFDE (\ref{sys}) is $p$th moment
ISS if there exist a continuous functional $V(t,\phi,i):J_{\tau}\times$
$\mathbf{C}\left(  \left[  -\tau,0\right]  ,\mathbf{R}^{n}\right)
\times \mathbf{S}\rightarrow J$, a USF $\mu(t)$, $\alpha_{1}\in \mathcal{VK}%
_{\infty},\alpha_{2}\in \mathcal{CK}_{\infty},$ and $\varpi \in \mathcal{K}%
_{\infty}$, such that the following conditions are met for all $t\in J$:
\begin{align*}
\mathrm{(a}_{1}\mathrm{).}\; \;  &  \alpha_{1}\left(  \left \vert
\phi(0)\right \vert ^{p}\right)  \leq V\left(  t,\phi,i\right)  \leq \alpha
_{2}(\left \Vert \phi \right \Vert ^{p}),\; \forall i\in \mathbf{S.}\\
\mathrm{(b}_{1}\mathrm{).}\; \;  &  \mathbf{E}\left \{  \mathcal{L}_{V}\left(
t,x_{t},i\right)  \right \}  \leq \mu(t)\mathbf{E}\left \{  V\left(
t,x_{t},i\right)  \right \}  ,\forall i\in \mathbf{S,}\quad \text{if}\;
\max_{i\in \mathbf{S}}\left \{  \mathbf{E}\left \{  V\left(  t,x_{t},i\right)
\right \}  \right \}  \geq \varpi \left(  \left \vert u\left(  t\right)
\right \vert \right)  \mathbf{.}%
\end{align*}

\end{theorem}

\begin{theorem}
\label{th4}(Krasovskii Theorem for iISS) The SFDE (\ref{sys}) is $p$th moment
iISS if there exist a continuous functional $V(t,\phi,i):J_{\tau}%
\times \mathbf{C}\left(  \left[  -\tau,0\right]  ,\mathbf{R}^{n}\right)
\times \mathbf{S}\rightarrow J$, a USF $\mu(t)$, $\alpha_{1}\in \mathcal{VK}%
_{\infty},\alpha_{2}\in \mathcal{CK}_{\infty},$ and $\varpi_{i}\in \mathcal{K}$,
$i=1,2$ such that $\mathrm{(a}_{1}\mathrm{)}$ and the following conditions are
met for all $t\in J$:
\[
\mathrm{(b}_{2}\mathrm{).}\; \; \mathbf{E}\left \{  \mathcal{L}_{V}\left(
t,x_{t},i\right)  \right \}  \leq \left(  \varpi_{1}\left(  \left \vert u\left(
t\right)  \right \vert \right)  +\mu(t)\right)  \mathbf{E}\left \{  V\left(
t,x_{t},i\right)  \right \}  +\varpi_{2}\left(  \left \vert u\left(  t\right)
\right \vert \right)  ,\; \forall i\in \mathbf{S}.
\]

\end{theorem}

Setting $u(t)\equiv0$ in Theorems \ref{th3} and \ref{th4} gives immediately
the following corollary regarding Krasovskii-type theorem for the $p$th moment
GUAS of the SFDE (\ref{sys}).

\begin{corollary}
\label{ppl}(Krasovskii Theorem for GUAS) The SFDE (\ref{sys}) is $p$th moment
GUAS if there exist a continuous functional $V\left(  t,\phi,i\right)
:J_{\tau}\times \mathbf{C}\left(  \left[  -\tau,0\right]  ,\mathbf{R}%
^{n}\right)  \times \mathbf{S}\rightarrow J\mathbf{,}$ a USF $\mu(t)$,
$\alpha_{1}\in \mathcal{VK}_{\infty}$ and $\alpha_{2}\in \mathcal{\mathbf{C}%
}\mathcal{K}_{\infty}$, such that $\mathrm{(a}_{1}\mathrm{)}$ and the
following conditions are met for all $t\in J$:
\[
\mathrm{(b}_{3}\mathrm{).}\  \mathbf{E}\left \{  \mathcal{L}_{V}\left(
t,x_{t},i\right)  \right \}  \leq \mu(t)\mathbf{E}\left \{  V\left(
t,x_{t},i\right)  \right \}  ,\; \forall i\in \mathbf{S}.
\]
If, in addition, Condition $\mathrm{(D)}$ is satisfied, then the system is
$p$th moment\ GUES.
\end{corollary}

Applying Theorems \ref{th3} and \ref{th4} and Corollary \ref{ppl} on system
(\ref{sys2}) gives the following corollary.

\begin{corollary}
\label{sys2Kras}Consider the SFDE (\ref{sys2}). Let $V\left(  t,\phi \right)
:J_{\tau}\times \mathbf{C}\left(  \left[  -\tau,0\right]  ,\mathbf{R}%
^{n}\right)  \rightarrow J$, $\mu(t)$ be a USF, $\alpha_{1}\in \mathcal{VK}%
_{\infty},\alpha_{2}\in \mathcal{CK}_{\infty},\varpi \in \mathcal{K},\varpi
_{i}\in \mathcal{K}_{\infty},i=1,2,q\in \mathcal{K}_{\infty}$, and $T\in
\Omega_{\mu}$. Consider the following conditions, where $t\in J$:
\begin{align*}
\mathrm{(a}_{2}\mathrm{).}  &  \; \; \alpha_{1}\left(  \left \vert
\phi(0)\right \vert ^{p}\right)  \leq V\left(  t,\phi \right)  \leq \alpha
_{2}(\left \Vert \phi \right \Vert ^{p}).\\
\mathrm{(b}_{4}\mathrm{).}  &  \; \; \mathbf{E}\left \{  \mathcal{L}_{V}\left(
t,x_{t}\right)  \right \}  \leq \mu(t)\mathbf{E}\left \{  V\left(  t,x_{t}%
\right)  \right \}  \mathbf{,}\quad \text{if}\; \mathbf{E}\left \{  V\left(
t,x_{t}\right)  \right \}  \geq \varpi \left(  \left \vert u\left(  t\right)
\right \vert \right)  .\\
\mathrm{(b}_{5}\mathrm{).}  &  \; \; \mathbf{E}\left \{  \mathcal{L}_{V}\left(
t,x_{t}\right)  \right \}  \leq \left(  \varpi_{1}\left(  \left \vert u\left(
t\right)  \right \vert \right)  +\mu(t)\right)  \mathbf{E}\left \{  V\left(
t,x_{t}\right)  \right \}  +\varpi_{2}\left(  \left \vert u\left(  t\right)
\right \vert \right)  .\\
\mathrm{(b}_{6}\mathrm{).}  &  \; \; \mathbf{E}\left \{  \mathcal{L}_{V}\left(
t,x_{t}\right)  \right \}  \leq \mu(t)\mathbf{E}\left \{  V\left(  t,x_{t}%
\right)  \right \}  .
\end{align*}
Then the SFDE (\ref{sys2})

\begin{enumerate}
\item is $p$th moment ISS if $\mathrm{(a}_{2}\mathrm{)}$ and $\mathrm{(b}%
_{4}\mathrm{)}$ are satisfied.

\item is $p$th moment iISS if $\mathrm{(a}_{2}\mathrm{)}$ and $\mathrm{(b}%
_{5}\mathrm{)}$ are satisfied.

\item is $p$th moment GUAS if $\mathrm{(a}_{2}\mathrm{)}$ and $\mathrm{(b}%
_{6}\mathrm{)}$ are satisfied.

\item is $p$th moment GUES if $\mathrm{(a}_{2}\mathrm{),(b}_{6}\mathrm{)}$ and
$\mathrm{(D)}$ are satisfied.
\end{enumerate}
\end{corollary}

By applying the Krasovskii functional based approach, the $p$th moment
asymptotic stability of SFDEs without Markovian switching and exponential
stability of stochastic delay interval systems with Markovian switching were
respectively investigated in \cite{Kolmanovskii1986} and \cite{Mao2002tac},
where time-derivatives of the Krasovskii functionals were required to be
negative definite. Compared with these results, in this section
time-derivatives of Krasovskii functionals can take both negative and positive
values along the solutions of SFDEs.

\section{\label{sec5}Stability of Deterministic Time-Delay Systems}

In this section, we discuss briefly the stability of the following
time-varying functional differential equation \cite{Ning2014tac}
\begin{equation}
\dot{x}\left(  t\right)  =f\left(  t,x_{t},u\left(  t\right)  \right)  ,\;t\in
J, \label{sys3}%
\end{equation}
where $f:J\times \mathbf{C}\left(  \left[  -\tau,0\right]  ,\mathbf{R}%
^{n}\right)  \times \mathbf{R}^{l}\rightarrow \mathbf{R}^{n}$ is assumed to be
locally Lipschitz in $\left(  t,x\right)  $ and uniformly continuous in $u,$
and to satisfy $f\left(  t,0,0\right)  =0,\forall t\in J.$ Clearly, this class
of systems is a special case of the SFDE (\ref{sys}) in the absence of
stochastic disturbance and Markovian switching. Similarly to the stochastic
setting, we give the following definitions \cite{Pepe2006SCL,Sontag2000tac}.

\begin{definition}
The trivial solution of system (\ref{sys3}) is said to be

\begin{enumerate}
\item ISS, if there exist $\sigma \in \mathcal{KL}$ and $\gamma_{1}%
\in \mathcal{K}$ such that, for any $u\in \mathbf{L}_{\infty}^{l},$
\[
\left \vert x\left(  t\right)  \right \vert \leq \sigma \left(  \left \Vert
x_{t_{0}}\right \Vert ,t-t_{0}\right)  +\gamma_{1}\left(  \left \Vert
u\right \Vert _{[t_{0},t]}\right)  ,\; t,t_{0}\in J,\;t\geq t_{0}.
\]

\item iISS, if there exist $\sigma \in \mathcal{KL}$ and $\gamma_{1},\gamma
_{2}\in \mathcal{K}$ such that, for any $u,$%
\[
\left \vert x\left(  t\right)  \right \vert \leq \sigma \left(  \left \Vert
x_{t_{0}}\right \Vert ,t-t_{0}\right)  +\gamma_{1}\left(  \int_{t_{0}}%
^{t}\gamma_{2}\left(  \left \vert u\left(  s\right)  \right \vert \right)
\mathrm{d}s\right)  ,\text{ }t,t_{0}\in J,\; t\geq t_{0}.
\]

\item GUAS if there exists a $\mathcal{KL}$-function $\sigma$ such that, when
$u\left(  t\right)  \equiv0,$
\[
\left \vert x\left(  t\right)  \right \vert \leq \sigma \left(  \left \Vert
x_{t_{0}}\right \Vert ,t-t_{0}\right)  ,\; t,t_{0}\in J,\;t\geq t_{0}.
\]

\item GUES if there exist two constants $\alpha>0$ and $\beta>0$ such that,
when $u\left(  t\right)  \equiv0,$%
\[
\left \vert x\left(  t\right)  \right \vert \leq \beta \mathrm{e}^{-\alpha \left(
t-t_{0}\right)  }\left \Vert x_{t_{0}}\right \Vert ,\;t,t_{0}\in J,\;t\geq
t_{0}.
\]

\end{enumerate}
\end{definition}

Applying Theorems \ref{th1} and \ref{th2} and Corollary \ref{GUAS} on system
(\ref{sys3}) gives the following Razumikhin-type stability theorem.

\begin{theorem}
\label{th5}(Razumikhin Theorem) Let $V(t,x):J_{\tau}\times \mathbf{R}%
^{n}\rightarrow J$ be a continuous function, $\mu(t)$ be a USF, $\alpha_{i}%
\in \mathcal{K}_{\infty},i=1,2,q\in \mathcal{K}_{\infty},\varpi \in
\mathcal{K}_{\infty}$, $T\in \Omega_{\mu}$, and $\rho \in(0,1)$ be a constant
satisfying (\ref{qs}). Consider the following conditions, where $t\in J$:
\begin{align*}
\mathrm{(H}_{1}\mathrm{).}\; \;  &  \alpha_{1}\left(  |x|\right)  \leq
V\left(  t,x\right)  \leq \alpha_{2}\left(  |x|\right)  .\\
\mathrm{(H}_{2}\mathrm{).}\; \;  &  \dot{V}\left(  t,x(t)\right)  \leq
\mu(t)V\left(  t,x(t)\right)  ,\\
&  \quad \text{if}\;V\left(  t,x(t)\right)  \geq \max \left \{  \varpi \left(
\left \vert u\left(  t\right)  \right \vert \right)  ,q^{-1}\left(  V\left(
t+\theta,x(t+\theta)\right)  \right)  \right \}  ,\; \forall \theta \in
\lbrack-\tau,0].\\
\mathrm{(H}_{3}\mathrm{).}\; \;  &  \dot{V}\left(  t,x(t)\right)  \leq
\mu(t)V\left(  t,x(t)\right)  +\varpi \left(  \left \vert u\left(  t\right)
\right \vert \right)  ,\quad \text{if}\;V\left(  t+\theta,x(t+\theta)\right)
\leq q\left(  V\left(  t,x(t)\right)  \right)  ,\; \forall \theta \in \left[
-\tau,0\right]  \mathbf{.}\\
\mathrm{(H}_{4}\mathrm{).}\; \;  &  \dot{V}\left(  t,x(t)\right)  \leq
\mu(t)V\left(  t,x(t)\right)  ,\; \text{if}\;V\left(  t+\theta,x(t+\theta
)\right)  \leq q\left(  V\left(  t,x(t)\right)  \right)  ,\; \forall \theta
\in \left[  -\tau,0\right]  \mathbf{.}%
\end{align*}
Then the time-varying time-delay system (\ref{sys3}) is

\begin{enumerate}
\item ISS, if $\mathrm{(H}_{1}\mathrm{)}$ and $\mathrm{(H}_{2}\mathrm{)}$ are satisfied.

\item iISS, if $\mathrm{(H}_{1}\mathrm{)}$ and $\mathrm{(H}_{3}\mathrm{)}$ are satisfied.

\item GUAS if $\mathrm{(H}_{1}\mathrm{)}$ and $\mathrm{(H}_{4}\mathrm{)}$ are satisfied.

\item GUES if $\mathrm{(H}_{1}\mathrm{)}$, $\mathrm{(H}_{4}\mathrm{)}$ and
$\mathrm{(D)}$ are satisfied.
\end{enumerate}
\end{theorem}

We mention that Item 3 of Theorem \ref{th5} is just Theorem 1 in
\cite{ZE16AUTO}, and Item 1 of this theorem improves Theorem 1 in
\cite{Ning2014tac} in the sense that $\mu(t)$ in Theorem \ref{th5} is less
restrictive than that in \cite{Ning2014tac} where the function needs to
satisfy condition (\ref{restrictive}).

Similarly, applying Theorems \ref{th3} and \ref{th4} and Corollary \ref{ppl}
on system (\ref{sys3}) gives the following theorem.

\begin{theorem}
\label{th6}(Krasovskii Theorem) Let $V(t,\phi):J_{\tau}\times \mathbf{C}\left(
\left[  -\tau,0\right]  ,\mathbf{R}^{n}\right)  \rightarrow J$ be a continuous
functional, $\mu(t)$ be a USF, $\alpha_{i}\in \mathcal{K}_{\infty}%
,i=1,2,\varpi \in \mathcal{K}_{\infty},$ and $\varpi_{i}\in \mathcal{K}$,
$i=1,2$. Consider the following conditions, where $t\in J$:
\begin{align*}
\mathrm{(h}_{1}\mathrm{).}\text{ \  \ }  &  \alpha_{1}\left(  \left \vert
\phi(0)\right \vert \right)  \leq V\left(  t,\phi \right)  \leq \alpha
_{2}(\left \Vert \phi \right \Vert ).\\
\mathrm{(h}_{2}\mathrm{).}\text{ \  \ }  &  \dot{V}\left(  t,x_{t}\right)
\leq \mu(t)V\left(  t,x_{t}\right)  ,\text{ if}\;V\left(  t,x_{t}\right)
\geq \varpi \left(  \left \vert u\left(  t\right)  \right \vert \right)
\mathbf{.}\\
\mathrm{(h}_{3}\mathrm{).}\text{ \  \ }  &  \dot{V}\left(  t,x_{t}\right)
\leq \left(  \varpi_{1}\left(  \left \vert u\left(  t\right)  \right \vert
\right)  +\mu(t)\right)  V\left(  t,x_{t}\right)  +\varpi_{2}\left(
\left \vert u\left(  t\right)  \right \vert \right)  .\\
\mathrm{(h}_{4}\mathrm{).}\text{ \  \ }  &  \dot{V}\left(  t,x_{t}\right)
\leq \mu(t)V\left(  t,x_{t}\right)  .
\end{align*}
Then the time-varying time-delay system (\ref{sys3}) is

\begin{enumerate}
\item ISS if $\mathrm{(h}_{1}\mathrm{)}$ and $\mathrm{(h}_{2}\mathrm{)}$ are satisfied.

\item iISS if $\mathrm{(h}_{1}\mathrm{)}$ and $\mathrm{(h}_{3}\mathrm{)}$ are satisfied.

\item GUAS $\mathrm{(h}_{1}\mathrm{)}$ and $\mathrm{(h}_{4}\mathrm{)}$ are satisfied.

\item GUAS $\mathrm{(h}_{1}\mathrm{)}$, $\mathrm{(h}_{4}\mathrm{)}$ and
$\mathrm{(D)}$\ are satisfied.
\end{enumerate}
\end{theorem}

Clearly, Item 3 of Theorem \ref{th6} is just Theorem 2 in \cite{ZE16AUTO}. We
point out that Items 1 and 2 of Theorem \ref{th6} generalize respectively
Theorems 1 and 3 in \cite{Ning2015AMC} in the sense that $\mu(t)$ here does
not need to satisfy condition (\ref{restrictive}) as required in
\cite{Ning2015AMC}.

\section{\label{sec6}Numerical Examples}

In this section we use two numerical examples to illustrate the obtained
theory. To same spaces, we only illustrate the Razumikhin-type stability theorems.

\textbf{Example 1}: Let $w\left(  t\right)  $ be a scalar Brownian motion and
$r\left(  t\right)  $ be a right-continuous Markov chain taking values in
$\mathbf{S}=\left \{  1,2\right \}  $ with generator%
\[
\Gamma=\left(  \gamma_{ij}\right)  _{2\times2}=%
\begin{bmatrix}
-1 & 1\\
2 & -2
\end{bmatrix}
.
\]
Assume that $w\left(  t\right)  $ and $r\left(  t\right)  $ are independent.
Let $d:J\times \mathbf{S}\rightarrow \lbrack0,\tau]$ be Borel measurable.
Consider the following system
\begin{equation}
\mathrm{d}x\left(  t\right)  =f\left(  t,x\left(  t\right)  ,r\left(
t\right)  \right)  \mathrm{d}t+g\left(  t,x\left(  t-d\left(  t,r(t)\right)
\right)  ,r\left(  t\right)  \right)  \mathrm{d}w\left(  t\right)  ,
\label{sysexample}%
\end{equation}
where $t\in J,$ and%
\begin{align*}
f\left(  t,x,1\right)   &  =-\frac{1}{2}x-\left \vert \sin t\right \vert
\sqrt[4]{\left \vert x\right \vert }\sqrt[3]{x},\;f\left(  t,x,2\right)
=-\frac{1}{2}x,\\
g\left(  t,y,1\right)   &  =-\sqrt{\left \vert b\left(  t\right)  \right \vert
}y\cos t,\;g\left(  t,y,2\right)  =\sqrt{\left \vert b\left(  t\right)
\right \vert }y\sin t,
\end{align*}
with $x=x\left(  t\right)  ,y=x\left(  t-d\left(  t\right)  \right)  $ and
$b\left(  t\right)  $ being a scalar function to be defined. Observe that%
\[
2xf\left(  t,x,i\right)  \leq-x^{2},\;g^{2}\left(  t,y,i\right)
\leq \left \vert b\left(  t\right)  \right \vert y^{2},\;i=1,2.
\]
To examine the stability of system (\ref{sysexample}), we construct a Lyapunov
function candidate $V\left(  t,x\left(  t\right)  ,i\right)  =x^{2}\left(
t\right)  ,i=1,2$ and calculate%
\begin{equation}
\mathcal{L}_{V}\left(  t,x_{t},i\right)  \leq-x^{2}+\left \vert b\left(
t\right)  \right \vert y^{2},\;i=1,2. \label{eqzb7}%
\end{equation}
If the classical Razumikhin theorem (Theorem 2 in \cite{Mao1999} and Theorem
3.1 in \cite{HuangMao2008tcsi}) is applied, it follows from (\ref{eqzb7}) that
the system is mean-square GUES if
\begin{equation}
\left \vert b\left(  t\right)  \right \vert <1,\quad t\in J. \label{classical R}%
\end{equation}

We next show how to improve this result by using Theorem \ref{GUAS}. Let
$q\left(  s\right)  =qs$ in condition (\ref{EVtq1}). Then, under the condition
$\min_{j\in \mathbf{S}}\left \{  \mathbf{E}\left \{  V\left(  t+\theta
,x(t+\theta),j\right)  \right \}  \right \}  \leq q\left(  \max_{i\in \mathbf{S}%
}\left \{  \mathbf{E}\left \{  V\left(  t,x(t),i\right)  \right \}  \right \}
\right)  ,\forall \theta \in \left[  -\tau,0\right]  \mathbf{,}$ we obtain from
(\ref{eqzb7}) that
\[
\mathbf{E}\left \{  \mathcal{L}_{V}\left(  t,x_{t},i\right)  \right \}  \leq
\mu(t)\mathbf{E}\left \{  V\left(  t,x(t),i\right)  \right \}  ,\;i=1,2,
\]
where $\mu \left(  t\right)  =-1+q\left \vert b\left(  t\right)  \right \vert .$
Now consider a periodic function $b\left(  t\right)  $ with period $\omega=1,$
and \cite{Mazenc2015CDC, ZE16AUTO}
\[
b\left(  t\right)  =\left \{
\begin{array}
[c]{c}%
0,\;t\in \lbrack0,c),\\
e,\;t\in \lbrack c,1),
\end{array}
\right.  \quad \left \vert b\left(  t\right)  \right \vert \leq e,
\]
where $c\in \left(  0,1\right)  $ and $e>0$ are some constant. We claim that
system (\ref{sysexample}) is mean-square GUES for arbitrary large $e$ if $1-c$
is sufficiently small.

It is easy to see that the function $\mu \left(  t\right)  $ is USF if and only
if \cite{ZE16AUTO}
\begin{equation}
\int_{0}^{1}\mu \left(  t\right)  \mathrm{d}t=-1+q\left(  1-c\right)  e<0.
\label{USFcon1}%
\end{equation}
If $e<1,$ then there exists a $q>1$ such that $\mu \left(  t\right)  <0.$
Consequently, $\mu \left(  t\right)  $ is USF for any $\tau$ and system
(\ref{sysexample}) is mean square GUES. Hence, without loss of generality, we
assume that $e\geq1.$ By (12) in \cite{ZE16AUTO}, we can compute%
\[
\varphi_{\mu}\left(  \omega \right)  =\int_{c}^{1}\mu \left(  t\right)
\mathrm{d}t=-1+q\left(  1-c\right)  e+c.
\]
Hence Condition (\ref{qs}) is equivalent to
\begin{equation}
q>\exp \left(  -1+q\left(  1-c\right)  e+c\right)  . \label{USFcon2}%
\end{equation}
These two inequalities (\ref{USFcon1}) and (\ref{USFcon2}) are equivalent to%
\begin{equation}
e<\frac{1}{\left(  1-c\right)  \exp \left(  c\right)  }. \label{USFcon3}%
\end{equation}
To summarize, system (\ref{sysexample}) is mean-square GUES if (\ref{USFcon3})
is satisfied.

Clearly, Condition (\ref{USFcon3}) is better than (\ref{classical R}) which is
equivalent to $e<1,$ in the sense that $e$ can be arbitrary large in the
former condition, by setting $1-c$ to be sufficiently small. Finally, we
mention that the results in the above cannot be obtained by the approach in
\cite{pz10tac} since the corresponding function $\mu \left(  t\right)  $ is
periodic and takes positive values in non-empty intervals.

\textbf{Example 2}: We consider the following time-varying SFDE%
\begin{equation}
\mathrm{d}x\left(  t\right)  =f\left(  t,x\left(  t\right)  ,r\left(
t\right)  ,u\left(  t\right)  \right)  \mathrm{d}t+g\left(  t,x\left(
t-d\left(  t,r(t)\right)  \right)  ,r\left(  t\right)  ,u\left(  t\right)
\right)  \mathrm{d}w\left(  t\right)  ,\text{ }t\in J, \label{examplesys3}%
\end{equation}
where $r\left(  t\right)  $ and $d(t,r(t))$ are the same as that in Example 1,
in which the generator $\Gamma$ is replaced by%
\[
\Gamma=\left(  \gamma_{ij}\right)  _{2\times2}=%
\begin{bmatrix}
-1 & 1\\
1 & -1
\end{bmatrix}
.
\]
The corresponding functions $f$ and $g$ are given by%
\begin{align*}
f\left(  t,x,1,u\right)   &  =-\frac{\lambda}{2}x-\frac{1}{2}x^{3}+\frac{t\cos
t^{2}}{1+x^{2}}u,\\
f\left(  t,x,2,u\right)   &  =\frac{1}{4}\left(  t\cos t^{2}-\lambda \right)
x+\frac{t\cos t^{2}}{2\left(  1+x^{2}\right)  }u,\\
g\left(  t,x,y,1,u\right)   &  =x^{2}-\sqrt{\frac{l}{2}}\sin^{k}\left(
t\right)  y,\\
g\left(  t,y,2,u\right)   &  =\frac{\sqrt{2l}}{2}\sin^{k}\left(  t\right)  y,
\end{align*}
with $x=x\left(  t\right)  ,y=x\left(  t-d\left(  t\right)  \right)  $ and
$u=u\left(  t\right)  $. Here $\lambda>\frac{1}{2\pi}$ and $l>0$ are
constants. Observe that%
\[
\left \{
\begin{array}
[c]{rl}%
2xf\left(  t,x,1\right)  & \leq-\lambda x^{2}-x^{4}+\frac{2t\cos t^{2}%
x}{1+x^{2}}u,\\
g^{2}\left(  t,y,1,u\right)  & \leq x^{4}+l\sin^{2k}\left(  t\right)  y^{2},
\end{array}
\right.
\]
and%
\[
\left \{
\begin{array}
[c]{rl}%
2xf\left(  t,x,1,u\right)  & \leq \frac{-\lambda}{2}x^{2}+\frac{t\cos t^{2}%
x}{1+x^{2}}u,\\
g^{2}\left(  t,y,1,u\right)  & \leq \frac{1}{2}l\sin^{2k}\left(  t\right)
y^{2}.
\end{array}
\right.
\]
To examine the stability of system (\ref{examplesys3}), we construct a
Lyapunov function candidate%
\[
V\left(  t,x\left(  t\right)  ,i\right)  =\left \{
\begin{array}
[c]{r}%
x^{2},\text{ }i=1,\\
\frac{1}{2}x^{2},i=2.
\end{array}
\right.
\]
Let $\left \vert u\left(  t\right)  \right \vert \leq V\left(  t,x\left(
t\right)  ,i\right)  $, $q\left(  s\right)  =qs$ and calculate%
\[
\mathcal{L}_{V}\left(  t,x_{t},i\right)  \leq \left(  t\cos t^{2}%
-\lambda \right)  x^{2}+ql\sin^{2k}\left(  t\right)  y^{2},\;i=1,2.
\]
Then, under Condition\textbf{ }$\mathrm{(B}_{1}\mathrm{)}$ in Theorem
\ref{th1}, we obtain from (\ref{ISSEV}) that
\[
\mathbf{E}\left \{  \mathcal{L}_{V}\left(  t,x_{t},i\right)  \right \}  \leq
\mu(t)\mathbf{E}\left \{  V\left(  t,x(t),i\right)  \right \}  ,\;i=1,2,
\]
where $\mu \left(  t\right)  =t\cos t^{2}-\lambda+ql\sin^{2k}\left(  t\right)
.$ According to Example 4.3 in \cite{ZE16AUTO}, if we set $q=\exp(2\lambda
\pi)$ and $k$ to be sufficiently large such that
\[
\frac{2k-1}{2k}\frac{2k-3}{2k-2}\cdots \frac{1}{2}<\frac{2\lambda \pi
-1}{2\lambda l\exp \left(  2\lambda \pi \right)  },
\]
then $\mu \left(  t\right)  $ is USF, $2\pi \in \Omega_{\mu}$, and $q>\exp
(\varphi_{\mu}\left(  2\pi \right)  ).$ Hence system (\ref{examplesys3}) is
mean-square ISS by Theorem \ref{th1}. Again, this result cannot be obtained by
the approach in \cite{pz10tac} since $\mu \left(  t\right)  $ cannot satisfy
the condition in (\ref{restrictive}). In fact, we have%
\begin{align*}
\int_{0}^{\infty}\max \{ \mu \left(  t\right)  ,0\} \mathrm{d}t  &  >\int
_{0}^{\infty}\max \{t\cos t^{2}-\lambda,0\} \mathrm{d}t\\
&  >\sum_{i=1}^{\infty}\int_{\sqrt{2i\pi-\frac{1}{4}\pi}}^{\sqrt{2i\pi
+\frac{1}{4}\pi}}\max \{t\cos t^{2}-\lambda,0\} \mathrm{d}t\\
&  \geq \sum_{i=1}^{\infty}\int_{\sqrt{2i\pi-\frac{1}{4}\pi}}^{\sqrt
{2i\pi+\frac{1}{4}\pi}}\max \left \{  \frac{\sqrt{2}}{2}t-\lambda,0\right \}
\mathrm{d}t\\
&  \geq \sum_{i=i^{\ast}}^{\infty}\int_{\sqrt{2i\pi-\frac{1}{4}\pi}}%
^{\sqrt{2i\pi+\frac{1}{4}\pi}}\left(  \frac{\sqrt{2}}{2}t-\lambda \right)
\mathrm{d}t\\
&  =\frac{1}{2}\pi \sum_{i=i^{\ast}}^{\infty}\left(  \frac{\sqrt{2}}{4}%
-\frac{\lambda}{\sqrt{2i\pi+\frac{1}{4}\pi}+\sqrt{2i\pi-\frac{1}{4}\pi}%
}\right) \\
&  >\frac{1}{4}\pi \sum_{i=i^{\ast}}^{\infty}\left(  \frac{\sqrt{2}}{2}%
-\frac{\lambda}{\sqrt{2i\pi-\frac{1}{4}\pi}}\right) \\
&  =\infty,
\end{align*}
where $i^{\ast}$ is the minimal integer such that $\frac{\sqrt{2}}{2}%
\sqrt{2i^{\ast}\pi-\frac{1}{4}\pi}-\lambda \geq0.$

\section{\label{sec7}Conclusion}

In this paper, generalized Razumikhin-type and Krasovskii-type stability
theorems on $p$th moment ISS, iISS, and GUAS were proposed for time-varying
stochastic time-delay systems with Markovian switching. Based on the general
It\^{o} formula and the stochastic analysis theory, some new Lyapunov
function(al)s based stability theorems were proposed by using properties of
uniformly stable functions and the improved comparison principle. The most
significant feature of the proposed results is that they allow
time-derivatives of the Razumikhin functions and Krasovskii functionals to
take both negative and positive values. The proposed results also improve the
related existing results on the same topic by removing some restrictive
conditions. Some improved stability criteria for deterministic time-delay
systems were also obtained.

\end{document}